\documentclass{elsart}

\usepackage{graphics}
\usepackage{epsfig}
\usepackage{amssymb}
\begin{document}
\begin{frontmatter}
\title{Model independent pre-processing of X-ray powder diffraction profiles}
\author{M. Ladisa\corauthref{cor2}}
\corauth[cor2]{Corresponding Author: Massimo Ladisa
(m.ladisa@ic.cnr.it)}
\ead{m.ladisa@ic.cnr.it} \ead[url]{www.ic.cnr.it}
\address{Istituto di Cristallografia (IC-CNR),
         Via Amendola 122/O,
         70126 Bari, Italy}
\author{A. Lamura}
\ead{a.lamura@ba.iac.cnr.it} \ead[url]{www.iac.cnr.it}
\address{Istituto Applicazioni Calcolo (IAC-CNR), Sezione di Bari,\\
          Via Amendola 122/D,
          70126 Bari, Italy}
\author{T. Laudadio}
\ead{Laudadio@esat.kuleuven.ac.be} \ead[url]{www.kuleuven.ac.be}
\address{Katholieke Universiteit Leuven,
         Department of Electrical Engineering, Division ESAT-SCD (SISTA),
         Kasteelpark Arenberg 10,
         3001 Leuven-Heverlee, Belgium}
\author{G. Nico}
\ead{g.nico@ba.iac.cnr.it} \ead[url]{www.iac.cnr.it}
\address{Istituto Applicazioni Calcolo (IAC-CNR), Sezione di Bari,\\
          Via Amendola 122/D,
          70126 Bari, Italy}

\begin{abstract}
Precise knowledge of X-ray diffraction profile shape is crucial in 
the investigation of the properties of matter in crystals powder. 
Line-broadening analysis is a pre-processing step in most of 
the full powder pattern fitting softwares. 
Final result of line-broadening analysis strongly depends on preliminary 
three steps: Noise filtering, removal of background signal and peak fitting.
In this work a new model independent procedure for two of the aforementioned 
steps (background suppression and peak fitting) is presented.
The former is dealt with by using morphological mathematics,
while the latter relies on the Hankel Lanczos Singular Value 
Decomposition technique.
Real X-ray powder diffraction (XRPD) intensity profiles of Ceria samples 
are used to test the performance of the proposed procedure. 
Results show the robustness of this approach and its capability of 
efficiently improving the disentangling of instrumental broadening. 
These features make the proposed approach an interesting
and user-friendly tool for the pre-processing of XRPD data.
\end{abstract}


\end{frontmatter}

\section{Introduction}

X-ray powder diffraction (XRPD) technique is nowadays 
a well known tool to study crystalline properties, which provide 
important information for applications in fields 
such as nanotechnology \cite{cervellino05a,cervellino05b}. 
All the applications benefit from a reliable pre-processing aimed at 
enhancing the quality of XRPD data. The pre-processing procedure consists
of four steps: Denoising, background suppression, peak fitting and 
signal deblurring, also known as line-broadening (see Fig. 1). 
Many techniques can be applied to remove noise from XRPD data 
(see ref. \cite{ladisa05} and references therein).
The step of background suppression is needed to emphasize the peak 
features of the sample. Traditional techniques are: Young's polynomials, 
Chebyshev approximation and linear interpolation \cite{jenkis96}.
Peak fitting is a very challenging step to extract information about 
the properties of polycrystalline powder \cite{estevez05}.
The final stage in XRPD data pre-processing is the signal deblurring aiming at 
disentangling the instrumental line-broadening out of data.
On this topic we spot a recent review in ref. \cite{balzar04}.
In this paper we propose a new approach for two of the above 
pre-processing steps of XRPD data. 
The signal background is removed by means of a procedure 
based on morphological mathematics. The peak shape profile fitting is 
carried out by using a subspace-based parameter estimation method called
Hankel Lanczos Singular Value 
Decomposition (HLSVD) technique \cite{laudadio02}. 
The main advantage of this approach is  
twofold: While both background suppression and peak shape profile fitting 
are model independent, the signal deblurring procedure, employed thereafter, 
crucially benefits from the model independence of the previous two steps. 
Indeed, our signal deblurring steps over the 
major drawback of most of the XRPD analysis tools available in the 
literature: The peak 
overlapping problem, namely the difficulty in singling the peak out 
of the full XRPD profile. This problem is mainly due to the model-induced bias 
in background/peak profile reconstruction. As to the denoising we use 
a wavelet based filter, a popular tool available in several software 
packages.
Real XRPD intensity profiles of Ceria samples \cite{balzar04} are 
used to test the performances. Four different raw datasets were used
in pairs \cite{data-ceria}. 
For each pair, one dataset was collected
on the annealed Ceria specimen (representing the instrumental broadening) 
and the other was collected on the broadened sample. The selected pairs are
those measured at the University of Birmingham (a high resolution X-ray
laboratory) and at the National Synchrotron Light Source (NSLS X3B1).

The structure of the paper is as follows. Section \ref{sec:method}
is devoted to the description of the proposed method and of the results.
Conclusions are drawn in Section \ref{sec:conclusions}.

\section{The method}
\label{sec:method}

Two different raw datasets were downloaded for the Ceria sample: 
The instrumental standard representing the instrumental broadening 
and the Ceria XRPD pattern of the broadened sample.

Our strategy in disentangling the profile broadening out of the 
experimental sample relies on a three step procedure which is 
sketched in the sequel.
\subsection{Denoising}
\label{sec:wavelets}
The noise was removed by applying wavelet transforms to the full 
XRPD spectrum and subtracted prior to the background suppression.
\par\noindent
Among the many applications of wavelets, signal denoising has been deeply 
investigated and the wavelet filter can be considered as the state of 
art on this subject. The discrete wavelet transform (DWT) is a linear 
operator which modifies the data vector in a similar way as the discrete 
Fourier transform (DFT). In both cases, the transform, given by a $N \times N$ 
matrix acting on the input $N$-vector data, is invertible \cite{press92}. The 
matrix entries are combinations of basis functions (the familiar 
sines and cosines in the case of DFT). 
Tables 1-2 summarize the main properties of some wavelet bases.
An interesting property of wavelet 
basis is the localization in both space and frequency domain. This means that 
they have a finite support or a decay in both domains. Regularity is another 
important property of wavelet basis. Regularity $r$ means that the $r^{th}$ 
derivative exists {\it almost everywhere} (see Table 1). 
In the case of Daubechies wavelet, regularity depends on the order $N$
(see Table 2). 
A wavelet is defined by particular set of numbers 
${\displaystyle \left\{ c_k \right\}_{k=0,\ldots,2N-1} }$, 
called wavelet filter coefficients.
They are determined by imposing the constraints of $N$ vanishing moments, 
${\displaystyle \sum_{k=0}^{2N-1}\ (-1)^k k^m c_k = 0}$ ($m=0,\ldots,2N-1$),
see Table 1, and orthogonality 
${\displaystyle \sum_{k=0}^{2N-1}\ c_k c_{k+2m} = 2 \delta_{0,m}}$. The 
coefficients $c_k$ characterize a low-pass filter while the coefficients 
$b_k = (-1)^k c_{2N-1-k}$
results in a high-pass filter. 
Coefficients $c_k$ give the entries of a $N \times N$ matrix, 
which 
iteratively applies to the $N$-data vector thus resulting
in the $N$-vector 
of detail coefficients. The whole procedure described above is the 
wavelet transform 
(see \cite{press92} for further details). 
Generally speaking, the denoising procedure involves three steps. 
The basic version of the procedure is the following: 
\begin{itemize}
\item Calculate the wavelet transform of XRPD profile and sort the components 
of the 
output vector by increasing frequency. This shall result in $N$-vector 
containing the 
XRPD profile average coefficient and a set of detail coefficients.
\item Noise thresholding, calculated on the highest frequency detail 
coefficient of the wavelet 
spectrum. 
\item Signal reconstruction by using the average coefficient and thresholded 
detail coefficients.
\end{itemize}
In this paper we choose the Daubechies wavelet basis with $N=2$. 
We address the reader to Daubechies \cite{daubechies92} for further details.
\subsection{Background suppression}
\label{sec:morphing}
The background was determined by means of morphological transforms 
for the full XRPD spectrum and subtracted prior to the signal deblurring.
\par\noindent
The morphological mathematics is based on the language of set theory. 
Considered a discrete binary image ${\mathcal I} \in \Re^2$ 
and a structuring element ${\mathcal S} \in \Re^2$, the four basic 
morphological mathematical operations on ${\mathcal I}$ by ${\mathcal S}$ 
are:
\begin{equation}
\begin{array}{lllll}
{\tt Dilation}: && {\mathcal I} \oplus {\mathcal S} &=& 
\cup_{s \in {\mathcal S}} {\mathcal I}_s \nonumber \\
{\tt Erosion}:  && {\mathcal I} \ominus {\mathcal S} &=& \cap_{s \in 
{\mathcal S}} {\mathcal I}_{-s} \nonumber \\
{\tt Opening}:  && {\mathcal I} \bigcirc {\mathcal S} &=& 
\left( {\mathcal I} \ominus{\mathcal S} \right) \oplus {\mathcal S} 
\nonumber \\
{\tt Closing}:  && {\mathcal I} \bullet {\mathcal S} &=& 
\left( {\mathcal I} \oplus {\mathcal S} \right) \ominus {\mathcal S} 
\nonumber \; 
\end{array}
\end{equation}
where ${\mathcal I}_s$ denotes the translation of ${\mathcal I}$ by $s$, 
namely 
${\displaystyle {\mathcal I}_s = \left\{ x+s | x \in {\mathcal I} \right\}}$. 
The value of each pixel in the output image is based on a comparison of the 
corresponding pixel in the input image with its neighbours, whose number and 
location 
is given by the structuring element. Generally, dilation expands 
image objects, whereas 
erosion shrinks them. In practice, dilation and erosion are employed in pairs. 
Opening is the erosion of an image followed by the dilation of the eroded 
image, and closing is the dilation of an image followed by the erosion of 
the dilated image. Opening eliminates sharp peaks smaller then the 
structuring element while the closing fills in the small holes and gaps. 
The binary morphological operations of dilation, erosion, opening and closing 
can be extended to grey-scale images. Let ${\mathcal E}_{\mathcal I}$ and 
 ${\mathcal E}_{\mathcal S}$ be the domains of the gray-scale image 
${\mathcal I}$ 
and the grey-scale structuring element ${\mathcal S}$ respectively. 
The grey-scale
dilation and erosion can be computed by
\begin{equation}
\begin{array}{lllll}
{\tt Dilation}: && \left( {\mathcal I} \oplus {\mathcal S} \right)(x,y) &=& 
\max \left\{ {\mathcal I}(x-m,y-n) + {\mathcal S}(m,n) \right\} \nonumber \\
{\tt Erosion}:  && \left( {\mathcal I} \ominus {\mathcal S} \right)(x,y) &=& 
\min \left\{ {\mathcal I}(x-m,y-n) - {\mathcal S}(m,n) \right\} \nonumber \;
\end{array}
\end{equation}
where $(x-m,y-n) \in {\mathcal E}_{\mathcal I}$ and   
$(m,n) \in {\mathcal E}_{\mathcal S}$. For such images, the minimum and 
maximum values are computed within neighbourhood  represented by the 
structuring 
element (see \cite{serra94} for details). 
\par
In our background suppression procedure the XRPD pattern is reshaped and 
padded into a 2-D image. A disk with a radius of three 
pixels is used as structuring element both for erosion and for  dilation. 
As to the erosion (dilation), pixels beyond the image border are assigned 
the maximum (minimum) value afforded by the data type. The morphological 
opening removes small objects from the image while preserving the shape 
and size of larger objects in the image. The overall result is a peak 
smearing effect while the background intensity remains unalterated. 
Restoring the original 1-D pattern provides the XRPD spectrum 
background. Figure 2 summarizes the whole procedure. 
We compared our findings to the traditional interpolation 
method and we found a satisfactory agreement (the percentage 
difference between the background computed by traditional techniques and 
our finding is below 3 \%). Up to our knowledge, 
this technique has never been applied to XRPD spectrum background 
suppression and it provides a reliable and user independent 
estimate of it.
\subsection{Peak fitting}
\label{sec:fitting}
The main problem in analysing an XRPD spectrum is the peak search, 
since the exact $\vartheta$ position is crucial in the extraction 
of the relevant microstructural information. Had the peak well defined, 
its shape would be straighforwardly achieved (for instance by a 
high resolution interpolation/fit by means of a model). Unfortunately 
the data resolution is rarely high enough to reach the goal of a 
well profiled peak. In that respect several methods have been devised so 
far to accomplish the peak fitting by means of gaussian, lorentzian, voigt, 
pseudo-voigt, Pearson VII and other models \cite{estevez05}. The main drawback 
of the aforementioned
methods cited above is the dependence of results on the model 
used in the fit procedure itself and the poor description of the real peak 
profile shape (for instance asymmetry). Here we use HLSVD 
method to the purpose \cite{laudadio02}. 
The main advantage of this 
method is the flexibility since the number of parameters 
is not fixed and it can be chosen to achieve a more satisfactory agreement 
between the model and the real peak profile shape. The HLSVD method works as
follows.
Let us model the XRPD intensity samples $I_n$ collected at angles 
$\vartheta_n$, $n = 0, \ldots N-1$ as the sum of $K$ exponentially damped 
complex sinusoids
\begin{equation}
I_n \simeq \sum_{k=1}^K a_k \exp (-d_k\ \vartheta_n) \ \cos[2 \pi f_k 
\vartheta_n + \varphi_k] \;,
\label{eq:model}
\end{equation}
$a_k$ is the amplitude, $\varphi_k$ the phase,
$d_k$ the damping factor and $f_k$ the frequency of the $k^{th}$
sinusoid, $k = 1, \ldots, K$, with $K$ the number of damped
sinusoids. 
The $N$ data points defined in (\ref{eq:model}) are arranged into
a Hankel matrix $H \stackrel{def.}{=} H_{L \times M}^{m\ l} = I_{m+l}$, 
$m=0,\ldots,M-1$, $l=0,\ldots,L-1$ , with $L+M = N+1$ 
($M \simeq L \simeq N/2$).
The SVD of the Hankel matrix  is computed as
$ H_{L \times M} = U_{L \times L} \Sigma_{L \times M} V^H_{M \times
M}$, where $\Sigma = {\rm diag}\{\lambda_1, \lambda_2, \ldots,
\lambda_r\}$, $\lambda_1 \geq \lambda_2 \geq \ldots \lambda_r$, $r
= \min(L, M)$, $U$ and $V$ are orthogonal matrices and the
superscript H denotes the Hermitian conjugate. The Lanczos 
bidiagonalization algorithm with partial
reorthogonalization %\cite{simon84} 
is used to compute SVD. 
This algorithm, based on FFT, computes the two matrix-vector
products which are performed at each step of the Lanczos procedure
in $O((L+M) log_2 (L+M))$ rather than in $O(LM).$
In order to obtain the ``signal'' subspace, the matrix $H$ is
truncated to a matrix $H_{K}$ of rank $K$
$H_K = U_K \Sigma_K V^H_K$, where $U_K$, $V_K$, and $\Sigma_K$ 
are defined by taking the first
$K$ columns of $U$ and $V$, and the $K \times K$ upper-left matrix
of $\Sigma$, respectively. As subsequent step, the least-squares
solution of the following over-determined  set of equations is
computed ${\displaystyle V^{(top)}_K E^H \simeq V^{(bottom)}_K}$, 
where ${\displaystyle V^{(bottom)}_K}$ and ${\displaystyle V^{(top)}_K}$ 
are derived from $V_K$ by
deleting its first and last row, respectively.
The $K$ eigenvalues ${\hat z}_k$ of matrix $E$ are used to
estimate the frequencies ${\hat f}_k$ and damping factors ${\hat
d}_k$ of the model damped sinusoids from the relationship
\begin{equation}
{\hat z}_k = \exp \left [ \left ( - {\hat d}_k + \imath 2 \pi
{\hat f}_k \right ) \Delta \vartheta \right ], \label{eq:eigen}
\end{equation}
with $k = 1, \ldots, K$. Values so obtained are inserted into the
model equation (\ref{eq:model}) which yields the set of equations
\begin{equation}
I_n \simeq \sum_{k=1}^K a_k \exp (-\hat d_k\ \vartheta_n) \ 
\cos[2 \pi \hat f_k \vartheta_n + \varphi_k] \;,
\label{eq:model_eigen}
\end{equation}
with $n = 0, \ldots, N-1$. The least-squares solution of
(\ref{eq:model_eigen}) provides the amplitude ${\hat a}_k$ and
phase ${\hat \varphi}_k$ estimates of the model sinusoids which are used 
in the next step.
\subsection{Deblurring}
\label{sec:deblurring}
The XRPD pattern has to be corrected for the instrumental 
broadening. Several methods have been devised so far to deal with 
this problem. Among them we quote the Stokes method \cite{stokes48} 
and the Bayesian approach \cite{lucy72,richardson72}. The main drawbacks 
of these methods stem from the difficulty in evaluating 
the background level, mainly due to peak overlapping. 
\par\noindent
The XRPD pattern plugged in the deblurring algorithm is noise-background 
free since it has been already pre-processed by the wavelets filter and the 
morphological operator.
\par\noindent
The technique proposed in this paper is a modified version of the 
one presented in ref. \cite{balzar04}. A blurred or degraded XRPD 
pattern can be approximately described by a Volterra equation 
$g = {\mathcal H} \otimes f + n$, 
where $g$ is the blurred XRPD pattern and ${\mathcal H}$ is the distortion 
operator due to several causes, namely the point spread function 
(PSF), and $n$ is an additive noise, introduced in the XRPD acquisition, 
that corrupts the signal. 
Strictly speaking in an XRPD experimental setup we deal with a poissonian 
noise which is a multiplicative noise. However the Poisson distribution 
function resembles the Gauss one provided a sufficiently large statistics 
in photons counting.
\par\noindent
As to the deblurring procedure we implement the damped Lucy--Richardson 
algorithm. This function performs multiple iterations, using optimization 
techniques and Poisson statistics. In our approach the PSF is the raw 
dataset downloaded for the Ceria sample - the instrumental standard - 
resembling the instrument profile \cite{balzar04}. 
The algorithm maximizes the likelihood that the resulting image, 
when convolved with the PSF, is an instance of the blurred image, assuming 
Poisson noise statistics. This function can be effective when the 
PSF is known but the knowledge about the additive noise in the image is poor. 
The Lucy--Richardson algorithm introduces several adaptations to the original 
maximum likelihood algorithm that addresses complex image restoration tasks. 
By using these adaptations, the effect of noise amplification 
on image restoration can be reduced, nonuniform image quality can be accounted 
for (e.g., bad pixels, flat-field variation) and the restored image resolution 
can be improved by subsampling.
\par\noindent 
Due to the denoising/background suppression, the original 
Volterra equation is readily simplified: $g={\mathcal H}\otimes f$, 
where $g$, ${\mathcal H}$ and $f$ have been already defined. Was the inverse 
${\mathcal H}^{-1}$ explicitly known, we would solve the former 
equation at a glance. Unfortunately this is not the case and the solution 
has to be approximated as follows. $g$, ${\mathcal H}$ and $f$ are positive 
and the PSF cannot change the $g$ norm, {\it i.e.} $||g||=||f||$. Thus:
\begin{equation}
\label{eq:lucy}
\sum_i f_i = \sum_i g_i 
\Longleftrightarrow \sum_i f_i = \sum_i g_i 
\frac{({\mathcal H}\otimes f)_i}{({\mathcal H}\otimes f)_i} 
= \sum_i f_i\ (\hat g \otimes {\mathcal H})^i \; ,
\end{equation} 
where ${\displaystyle \hat g_k = \frac{g_k}{({\mathcal H}\otimes f)_k}}$; 
the solution can be found by an iterative procedure: 
$\displaystyle{f^{n+1}=f^n\ \frac{g \otimes 
{\mathcal H}}{{\mathcal H}\otimes f^n}}$, 
where the initial guess XRPD spectrum is uniform. 
As already stressed in ref. \cite{balzar04}, the main 
drawbacks in applying such an algorithm to the single peak deconvolution are 
the noise amplification and the peak fitting bias. Noise amplification is 
dramatically reduced by both the denoising procedure and the small (some 
five) number of iterations used in the algorithm. Moreover the damp in the 
algorithm specifies the threshold level for the deviation of the resulting 
image from the original image, below which damping occurs. For pixels that 
deviate in the vicinity of their original values, iterations are suppressed. 
\par\noindent
As to the peak fitting bias, unlike the Balzar approach, our procedure uses 
the instrumental standard pattern (with no overlapping) as the PSF to 
deconvolve the full XRPD pattern and then we extract the 
deconvoluted/deblurred XRPD pattern in the same range of the PSF used for 
the deconvolution itself. The rationale of this choice relies on the fact 
that while the PSF peaks have no overlap, this is not the case for the 
broadened sample peaks and, thus, the peak ranges can be defined starting 
on the annealed sample rather than the broadened one.  Moreover the discrete 
Fourier transform (DFT), used by the deblurring functions, assumes that the 
frequency pattern of an image is periodic. This assumption creates a 
high-frequency drop-off at the edges of an overlapping peaks cluster 
\cite{biemond90}. 
This high-frequency drop-off can create an effect called boundary related 
ringing 
in deblurred images, that is a systematic error affecting any further 
investigation on the physical meaning of the deconvolved spectrum. To reduce 
ringing, our full pattern deconvolution, as described above, resembles an 
edgetaper function. It removes the high-frequency drop-off at the edge of an 
image by blurring the entire image and then replacing the center pixels of 
the blurred image with the original one. In this way, the edges of the 
image taper off to a lower frequency.

\section{Conclusions}

\label{sec:conclusions}
In this paper we presented a new approach for background removal and peak 
fitting
of XRPD profiles. Such operations are crucial in the line-broadening analysis 
of X-ray 
diffraction profiles, an important pre-processing step in the investigation 
of the crystal powder samples by means of XRPD data. 
Backgroud suppression relies on the use of morphological  mathematics while 
peak fitting is carried out by means of HLSVD technique.
In order to enhance the signal-to-noise ratio of XRPD profiles a wavelet 
based filter is 
preliminarly applied to XRPD data. 
The output of the proposed precedure is supplied to a damped Lucy--Richardson 
algorithm 
for deblurring.
The main advantage of this approach is twofold: Background suppression and 
peak shape 
profile fitting are model independent.
Real XRPD intensity profiles of Ceria samples are used to test performances. 
Results 
show that the background estimate is in agreement with that computed by 
traditional 
interpolation methods with a percentage difference below 3\%.
Further, the output XRPD profile has narrower peaks, located at the same 
position and 
with the same shape as the original one. 
It is worth noting that once the deblurred, noise-background free 
XRPD spectrum is convoluted back to the PSF and added to the removed noise 
and background signals, the resulting XRPD pattern resembles the original one.

\newpage
\begin{center}
\begin{table}[ht!]
\label{tab:1}
\caption{Properties of some wavelet bases.}
\begin{tabular}{l||c|c|c|c} 
 wavelet        & t-localization & f-localization &  \# zero moments &  $r$       \\
\hline \hline
 Haar           &      [0,1]     &       1/f      &        1         &  0         \\ 
 Sinc           &      1/t       &      [0,1]     &     $\infty$     & $\infty$   \\ 
 Daubechies (N) &      [0,2N-1]  &       1/f      &       N          & $\alpha$(N)\\
\hline
\end{tabular}
\end{table}
\end{center}
\begin{center}
\begin{table}[ht!]
\label{tab:2}
\caption{Regularity $r = \alpha (N)$ of a Daubechies wavelet basis of order $N$.}
\begin{tabular}{l||c|c|c|c|c|c} 
 N           &   2   &   3   &   4   &   5   &   6   &     N     \\
\hline \hline
 $\alpha$(N) & 0.500 & 0.915 & 1.275 & 1.596 & 1.888 & 0.2075\ N \\
\hline
\end{tabular}
\end{table}
\end{center}

\newpage
\begin{center}
\begin{figure}
\caption{Block diagram of the pre-processing procedure. The proposed
approach refers to the background suppression and peak fitting steps 
(shadowed boxes).}
\includegraphics*[width=\textwidth]{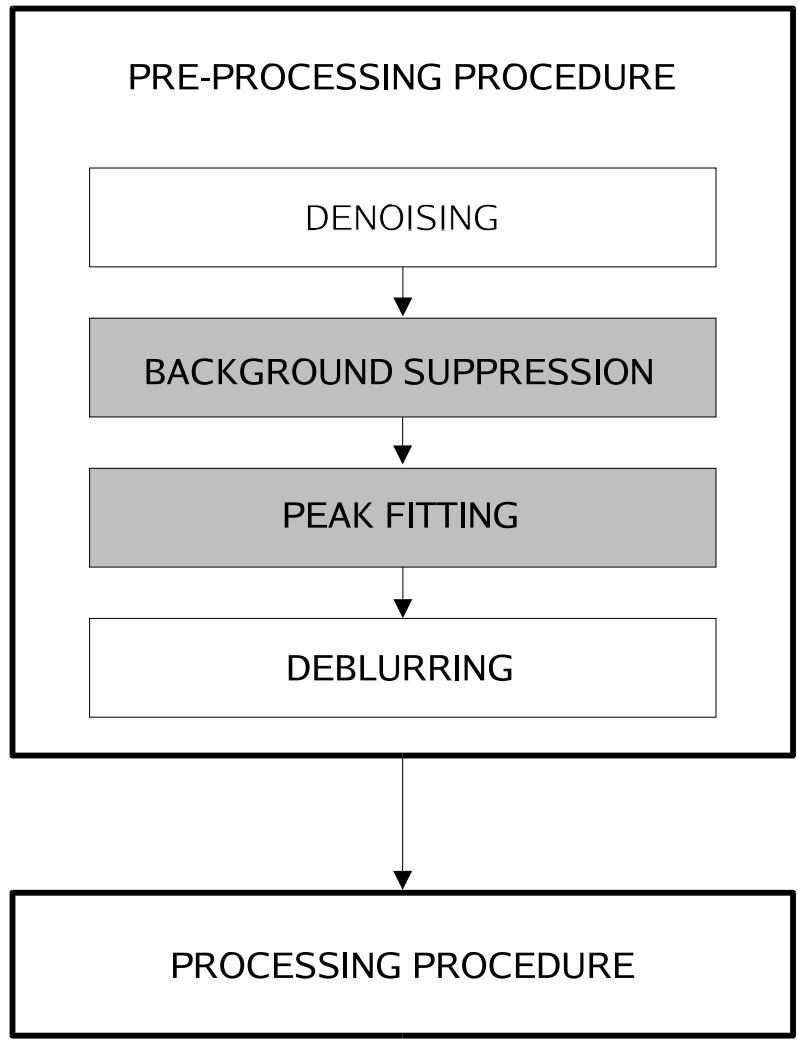}\\*
\end{figure}
\end{center}

\begin{center}
\begin{figure}
\caption{Background evalutation procedure. 
Left top: original XRPD pattern. 
Right top: XRPD pattern after two-dimensional reshaping.
Right bottom: two-dimensional reshaped XRPD pattern after morphological 
operations.
Left bottom: background after one-dimensional reshaping.}
\includegraphics*[width=\textwidth]{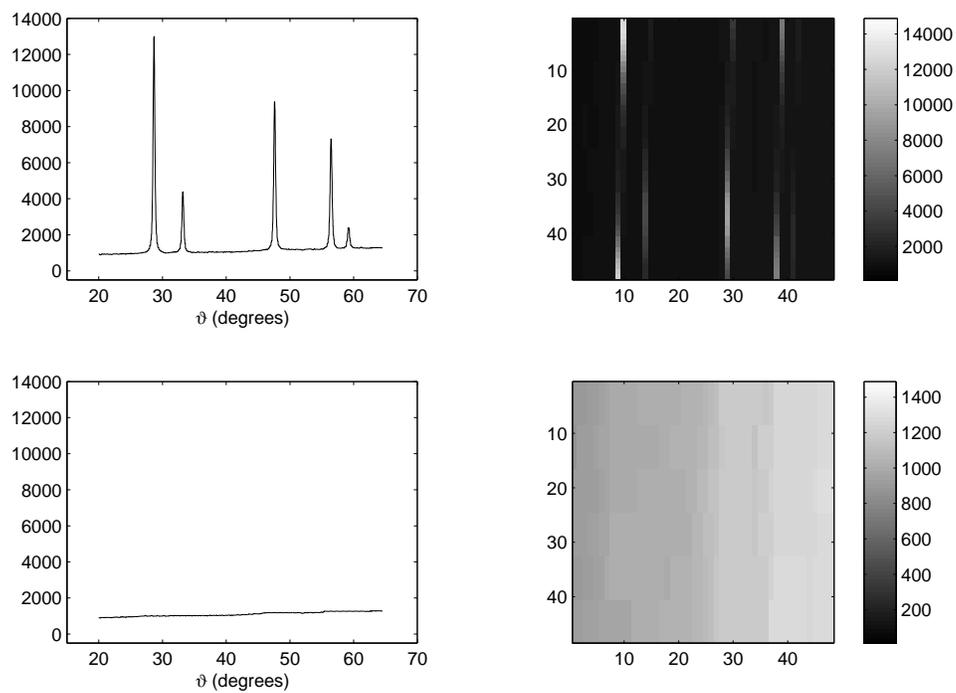}\\*
\end{figure}
\end{center}
\begin{center}
\begin{figure}
\caption{Plots. 
Left top: original XRPD pattern. 
Right top: XRPD pattern after denoising and background suppression.
Left bottom: noise-background free XRPD pattern after deblurring.
Right bottom: residue between the final XRPD pattern re-convoluted to the PSF 
together with noise and background and the original XRPD pattern.}
\includegraphics*[width=\textwidth]{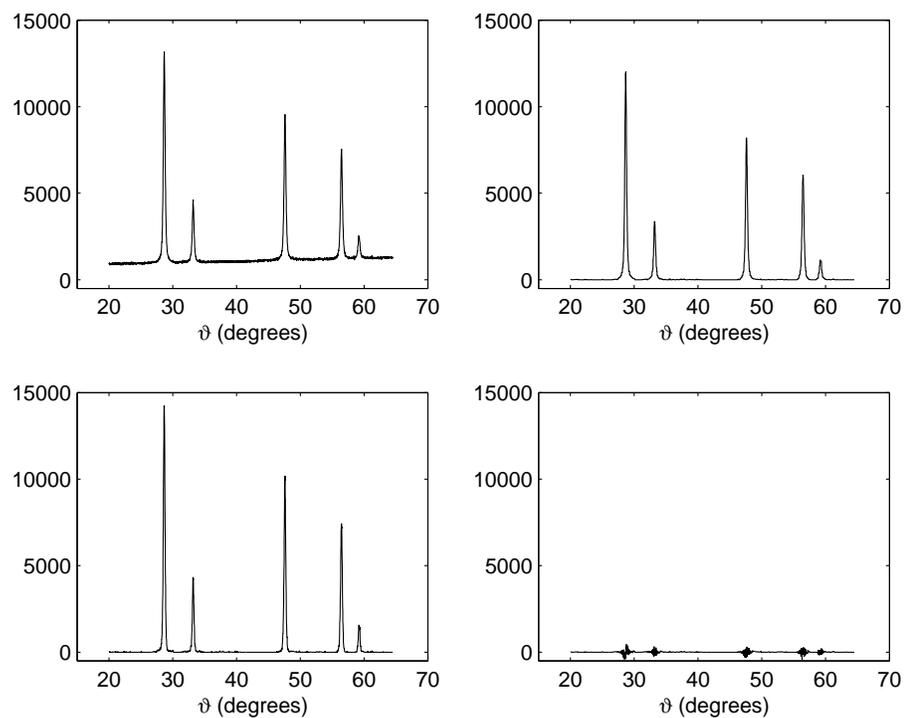}\\*
\end{figure}
\end{center}

\end{document}